\documentstyle[12pt,leqno]{article}
\title{On discriminants of polylinear forms}
\author{Valeri V.Dolotin}
\date{}
\begin{document}
\newtheorem{example}{Example}[subsection]
\newtheorem{const}{Construction}
\newtheorem{Definition}{Definition}[subsection]
\newtheorem{Proposition}{Proposition}[subsection]
\newtheorem{Theorem}{Theorem}[subsection]
\newcommand{\C}{\mbox{${\bf C}$}}
\newcommand{\PP}{\mbox{${\bf P}$}}
\def\theequation{\arabic{section}.\arabic{equation}}
\maketitle
\begin{abstract}
In this paper we propose a conseptual framework for the observed
properties of discriminants of polylinear forms. The connection with
classical problems of linear algebra is shown. A new class of
algebraic varieties (hypergrassmanians) is introduced, the particular case of
which are grassman manifolds. An algorithm is given for computing the
discriminants of polylinear forms of "boundary format" (hypergrassmanian
analogues of plukker coordinates). The algorithm for computing the
discriminants of polylinear forms of general formats is outlined.
\end{abstract}

\setcounter{section}{-1}
\section{Introduction}

It is convenient (and instructive) to keep for certain objects of
multilinear algebra the terminology from their counterparts in linear
algebra (see [1]). So for the set of elements of $d-$linear form will
be used the term "$d-$dimensional matrix", and for its discriminant
(see definition below) we can use the term "determinant" with the notion
of "minors" also making sense.

The investigation of determinants of multidimensional
matrices may be useful already in the linear algebra (see Section 1), which
essentially deals with sets of rectangular matrices. A good example
(Section 1.1) is
a 3-dimensional interpretation of the theory of kronekker pairs, which
in this context obtains a straightforward generalization. It is
also interesting to look for a generalization of eigenvalue theory.
The theory of eigenvalues in its different versions is
equivalent to the investigation of  matrices of type $A+\lambda B$,
where $A$ and $B$ is a paire of $n\times n$ matrices,  and
invariants of this pair, which can be expressed in terms of
$GL(n)\times GL(n)$-action on $n\times n\times 2$ form (3-dimensional
matrix). In a similar way the
``multidimansional eigenvalue theory'' is reduced to the invariants of
3-dimensional matrices of larger formats.

We know, that in the case of bilinear forms (2-dimansional matrices)
the determinant (as one polynomial) is defined only for square matrices.
So in order to specify the size of 2-dimensional matrix which has a
unique expression (determinant) as a characteristic of its degeneracy,
it is enough to give one number $n$ - the number of rows or columns.
For polylinear forms the difference from bilinear ones we can
illustrate on the case of 3 dimensions.
Take a 3-dimensional matrix with elements
$a_{i_1i_2i_3}$, where $1\le i_1\le n_1$, $1\le i_2\le n_2$, $1\le i_3\le n_3$.
Let us fix $n_1$ and $n_2$. Then the determinant is defined for matrices
with $n_3$ satisfiing the following inequality (see [1]):
%(1)
\begin{equation}
n_1-n_2+1\le n_3\le n_1+n_2-1
\end{equation}
So the size ("format") of $d$-dimensional matrices which have the determinant
is defined in general by $d$ parameters while in 2-dimentional case by only
1. But in d-dimensianl case for $d>2$ there also exists a class of
matrices the size of which is described by $d-1$ parameters. These are so
called ``matrices of boundary format'', the size of which corresponds
to the equality in $(0.1)$. In 2-dimensional case the matrices, the
degeneracy of which can not be characterized by one expression, are
rectangular matrices. In higher dimensional case these are ``matrices
of grassman format'', those which in case $d=3$ do not satisfy $(0.1)$.
It this paper we give a set of properties of matrices of boundary and
grassman formats, which show that they are the proper
generalizations if square and rectangular matrices correspondingly.
In particular, the condition for rectangular matrices of beeing of
corank 1 is that the determinants of all maximal square submatrices
(maximal minors) are 0. The corresponding first degeneracy condition for
d-dimensioanl matrices, defined here as the condition of beeing of
corank 1 (see Section 2.1), is that the determinants of all maximal
submatrices of boundary format are 0. To a rectangular matrix one can
put into correspondense a set of vectors - its rows or columns.
Then the condition for
beeing of corank 1 for rectangular matrix has a geometric
interpretation as the linear dependence of these vectors
(1-dimensional matrices). To a d-dimensional matrix with elements
$a_{i_1...i_d}$ one can put into correspondence a set of "slices" in
k-th direction, which are $(d-1)$-dimensional matrices with elements
$a_{i_1...\widehat i_k...i_d}$. Then the condition of
beeing of corank 1 for d-dimensional matrices of grassman format can be
expressed geometrically in terms of singularities of the intersection
of the span of these slices with the submanifold of
$(d-1)-$dimensional matrices of corank 1. The remarkable fact, making
the notion of corank 1 matrices well defined, is that this singularity
condition does not depend on the direction of slicing of our
matrix (the number k above).

Consider the problem of finding the kernel of a lenear combination
$S(\lambda):=\lambda_1 A_1+...+\lambda_k A_k$ of
k rectangular matrices of size $m\times n$. This kernel will be an
$(n-m)$-dimensional subspace, i.e. an element of $G_{m,n}$. Changing
$(\lambda_1,...,\lambda_k)$ we get an k-parametric subset in
$G_{m,n}$. In the case when $k=n-m+1$ the image of this subset via
plukker embedding of $G_{m,n}$ will be a Veronese manifold. The
Veronese manifolds which can be obtained in this way are called here
``proper Veronese manifolds''. The condition for a proper Veronese manifold to
be degenerate can be expressed in to ways:\\
$1)$ as the condition that there exists such $(\lambda_1,...,\lambda_k)$
that all plukker coordinates of the kernel of $S(\lambda)$ (which
are $m\times m$ minors of $S(\lambda)$) are equal to 0, or as the
singularity of the intersection of $span(A_1,...,A_k)$ with the
submanifold $M^{'}_{mn}$ of degenerate $m\times n$ matrices\\
$2)$ as the condition that the determinant of 3-dimensional matrix of
size $m\times n\times (n-m+1)$ made of the elements of
$A_1,...,A_{n-m+1}$ is 0.

In the case of grassman format, when $k>n-m+1$, the condition of the
existence of the intersection of $span(A_1,...,A_k)$ with $M^{'}_{mn}$
is that the determinants of all $m\times n\times (n-m+1)$ submatrices
of the corresponding 3-dimensional matrix are equal to 0. It happens
that the similar fact holds in general d-dimensional case. This allows
to interprete the determinants of maximal submatrices of boundary
format of matrix of grassman format as multidimensional anologues
of plukker coordinates and to consider the anologue of plukker map on
the space of d-dimensional matrices of grassman format:
$M_{n_1...n_d}\to \PP ^{\left( \matrix{n\cr m}\right)}$, putting into
correspondence for a matrix the set of its maximal minors of boundary
format. As in 2-dimensional case here arises the fundamental problem
to find the relations between the minors of multidimensional matrix,
the analogues of plukker relations, i.e. to discribe the image of
$M_{n_1...n_d}$ as an algebraic manifold.

For studiing the discriminants of polylinear forms there is a fundamental
question about the algorithm  of explicit calculation of these discriminants.
In Section 3 we develop a technique which gives an algorithm of calculating
the discriminants of $d-$linear forms of boundary format (``hyperplukker
koordinates''). This technique happens to be basic for calculating
the discriminants of $d-$linear forms of general format.
In Section 4 we outline and give an example of this general algorithm.

\begin{Definition}
{\rm If $p(x_1,...,x_m)=\displaystyle
\sum^{}_{i_1\le{...}\le{i_n}}c_{i_1{...}{i_{n}}}x_{i_{1}}
{...}x_{i_{n}}$ is a homogeneous polynomial of degree $n$,
then the set of values of coefficients $c$ is called {\em discriminantal}
if the system of equations
%(2)
\begin{equation}
\frac{\partial{p(x_1,...,x_m)}}{\partial{x_i}}=0,\quad  i=1,...,m
\end{equation}
has a solution $(x_1,...,x_m)\in{(\C^{*})^m}$.\\
In case when the discriminantal set is an algebraic submanifold of codimension
$1$ in the space of coefficients, it is called {\em the discriminant} of $p$,
denoted by $D(p)$.}
\end{Definition}
Now we will consider a particular case of this definition.
Let $V_{n_1},...,V_{n_d}$ be a set of linear vector spaces, such that
$dim(V_{n_i})=n_i$.
Let  $a\in (V_{n_1}\otimes ...\otimes V_{n_d})^*$ be a $d-$linear
form. For a set of vectors $x^{(k)}\in V_{n_k}$, $k=1,...,d$ with coordinates
$x^{(k)}=(x^{(k)}_1,...,x^{(k)}_{n_k})$ in a chosen basis the value of
the form on $\otimes^d_{k=1}x^{(k)}$ is a polinomial of $d$ sets of
variables $x^{(1)},...,x^{(d)}$ of degree $d$
$$a(x^{(1)},...,x^{(d)})=p(x^{(1)},...,x^{(d)})$$
and the system $(0.2)$ becomes
%(3)
\begin{equation}
\frac{\partial{a (x^{(1)},...,x^{(d)})}}{\partial{x^{(k)}_i}}=0,
\quad  i=1,...,n_k, \quad k=1,...,d
\end{equation}
The coefficients $a_{i_1...i_d}$ of this form are the elements of
$d-$dimensional rectangular $n_1\times ...\times n_d$ matrix.

\begin{Definition}
{\rm The discriminant of the polinomial $a(x^{(1)},...,x^{(d)})$ is
called {\em the determinant} of $n_1\times ...\times n_d$ matrix
$(a_{i_1...i_d})$ and denoted by $det(a)$.}
\end{Definition}
{\bf Example.} Let $a\in (U_{n}\otimes V_{n})^*$. Then its
coefficients form a usual $n\times n$ square matrix. Then
$h=\displaystyle\sum^{}_{1\le i,j\le n} x_i a_{ij} y_j$.
The system $(0.2)$ in this case becomes:
%(4)
\begin{equation}
\sum_{j=1}^n a_{ij} y_j=0,\quad i=1,...,n\quad
\sum_{i=1}^n x_ia_{ij}=0,\quad j=1,...,n
\end{equation}
Note, that this system contains the system of linear homogeneous
equations as well as its conjugate. For the discriminants of
polylinear forms this property will become essential.
There is a natural action of the group $GL_{n_1}\times ...\times GL_{n_d}$
on $\otimes^d_{j=1}V_{n_j}$ with the induced action on
$(\otimes^d_{j=1}V_{n_j})^*$. Since the system $(0.2)$ for
$p=a(x^{(1)},...,x^{(d)})$ is invariant under this action we have
\begin{Proposition}\label{P1}
The determinant of $n_1\times ...\times n_d$ matrix is invariant under
the action of $GL_{n_1}\times ...\times GL_{n_d}$.
\end{Proposition}
{\bf Notation} Denote $M_{n_1...n_d}:=\displaystyle
(\otimes^d_{j=1}V_{n_j})^{*}$.

%%%   SECTION 1
\setcounter{section}{0}
\section{Problems of linear algebra and determinants of
$3-$dimensional matricies}
\subsection{General setting}
Let $M_{nm}$,where $n\le m$, be a linear space of $n\times
m$ matricies. Let $M^{\prime}_{nm}\subset M_{nm}$ be a
submanifold of matricies of rank $n-1$. Let $A_1,...,A_k\in M_{nm}$
be a set of $n\times m$ matrices. From their elements we can
make an $n\times m\times k$ matrix of coefficients of 3-linear
form $(a_{i_{1}i_{2}i_3})$. There is $1-1$ correspondence between
linear subspaces $span(A_1,...,A_k)\subset M_{nm}$ for different
choices of $A_1,...,A_k$ and the orbits
of the corresponding forms $a (A_1,...,A_k)$ under
the action of $GL_k$ on $M_{nmk}=(V_n\otimes V_m\otimes V_k)^*$.

\begin{example}
Let $A$ and $B$ be $n\times n$ matrices.
\begin{Proposition}\label{P1}
Let $det(B)\ne 0$. Then the following statements are aquivalent:\\
1) $det(a_{i_{1}i_{2}i_{3}}(A,B))=0$\\
2) $D(det(A+zB))=0$ or
$D(det(AB^{-1}-zI))=0$, i.e. the  characteristic polinomial of
$AB^{-1}$ has multiple roots\\
where $D(p(z))$ denotes the discriminant of polinomial $p(z)$.
\end{Proposition}
\end{example}

\begin{example}
Let $A$ and $B$ be $n\times (n+1)$ matrices.
\begin{Proposition}
The determinant of $n\times (n+1)\times 2$ matrix
$det(a_{i_{1}i_{2}i_{3}}(A,B))\ne 0$ iff the pair $(A,B)$ is
``kronekker'', i.e. by the action of $GL_n\times GL_{n+1}$ it can be
reduced to the following canonical form\\
$$A=\left( \matrix{
1 & \      &{0} &0     \cr
\ & \ddots & \  &\vdots\cr
0 & \      &1   &0
}\right),\quad
B=\left( \matrix{
0      & 1  & \      & 0 \cr
\vdots & \  & \ddots & \  \cr
0      & 0  & \      & 1
}\right)$$
\end{Proposition}
\end{example}

\subsection{Veronese manifolds}
Let $k=m-n+1$, $A_1,...,A_k\in M_{nm}$.
Let $S(z_1,...,z_k):=z_1 A_1+...+z_k A_k$ be a point of subspace
$span(A_1,...,A_k)\subset M_{nm}$.
Let $\Delta(z_1,...,z_k):=(\Delta_{i_{1}...i_{m}}(S))_{(i_1,...,i_k)\in
\left( \matrix{[m]\cr n}\right)}\in \PP ^{\left( \matrix{n\cr
m}\right)}$ be a vector with components - $n\times n$ minors of $S$,
the image of $S(z)$ via Plukker embedding.
The map $\varphi : (z_1,...,z_k)\mapsto \Delta(z_1,...,z_k)$ gives a
$k-$parametric submanifold ${\cal{V}}(A_1,...,A_k)$ in
$\PP ^{\left( \matrix{m\cr n}\right)}$.

\begin{Proposition}
For $A_1,...,A_k$ in general position the intersection
$M^{\prime}_{nm}\cap span(A_1,...,A_k)$ is empty.
\end{Proposition}

\begin{Proposition}
For $A_1,...,A_k$ in general position the manifold
${\cal{V}}(A_1,...,A_k):=\varphi (\C^k)\in \PP ^{\left(
\matrix{m\cr n}\right)}$ is a veronese manifold.
\end{Proposition}
The veronese manifolds obtained in this way we call {\em proper}.

\begin{Theorem}
The following statements are equivalent:\\
1) the intersection $M^{\prime}_{nm}\cap span(A_1,...,A_k)$ is not empty\\
2) $det(a_{i_{1}i_{2}i_{3}}(A_1,...,A_k))=0$\\
3) the manifold ${\cal{V}}(A_1,...,A_k)$ is singular\\
4) ${\cal{V}}(A_1,...,A_k)$ belongs to a hyperplane in $\PP ^{\left(
\matrix{m\cr n}\right)}$
\end{Theorem}

Let $(x_1,...,x_n),(y_1,...,y_m),(z_1,...,z_k)$ be the coordinates in
$V_n,V_m$ and $V_k$ correspondingly.
For $a=a(A_1,...,A_k)\in M_{nmk}=V_n\otimes V_m\otimes V_k$ the
system $(0.2)$, where $p=a (x,y,z)$, contains a subsystem
%(1.5)
\begin{equation}
\frac{\partial{a(x,y,z)}}{\partial{x_{i}}}=
\sum^m_{j=1}(z_1A_1+...+z_kA_k)_{ij}y_j=S(z)y=0,\quad i=1,...,n
\end{equation}

To a given $z\in (\C^*)^k$ we can put into correspondence a subspace
$Ker(S(z))\subset V_m$ of solutions of (1.5). If $rank(S(z))=n$ then
$Ker(S(z))$ has dimension $(m-n)$. If $rank(S(z))<n$
then $Ker(S(z))$ has dimension greate then $(m-n)$.
\begin{Proposition}
$det(a(A_1,...,A_k))=0$ iff there are values of $z\in (\C^*)^k$
such that the dimension of the space of solutions of corresponding
system (1.5) is greate then $(m-n)$ and the $\varphi$-
image of the set of such $z$ is exactly the set of singular points of
${\cal{V}}(A_1,...,A_k)$.
\end{Proposition}

%%%%%%%%%%%%%%%%%%%%%%%%%%   SECTION 2

\section{Hyperveronese and hypergrassmanian manifolds}
\subsection{On the rank of polylinear forms}
For a given form $a\in M_{n_1...n_d}$ and an integer number $1\le k\le
d$ we can put into correspondence a set of forms ${a}^{(k)}_{i_k}\in
M_{n_1...\widehat{n}_k...n_d},\quad i_k=1,...,n_k$, such that
$({a}^{(k)}_{i_k})_{i_{1}...\widehat{i}_{k}...i_{d}}=
a_{i_{1}...i_{d}}$.
For a given $k$ there is $1-1$ correspondence between the orbits
of the forms $a$ under
the action of $GL_{n_k}$ on $M_{n_1...n_d}$
and the corresponding linear subspaces
$span({a}^{(k)}_1,...,{a}^{(k)}_{n_k})\subset M_{n_1...\widehat{n}_k...n_d}$.
We say, that the {\em $GL_{n_k}$-orbit of $a$ has an intersection} with
a submanifold in $M_{n_1...\widehat{n}_k...n_d}$ if
$span({a}^{(k)}_1,...,{a}^{(k)}_{n_k})$ has an intersection
with this submanifold.

{\bf Notation.} Denote by $M'_{n_1n_2}$ the set of $n_1\times n_2$ matrices of
corank 1.

For a multiindex $n_1...n_d$ let $\widehat{n}_k:=n_1...n_{k-1}n_{k+1}...n_d$.
Now by induction on $d$ we can introduce the following definition.

\begin{Definition}
A subset of the space $M_{n_1...n_d}$ is called the set
of forms of corank $1$ (denoted by $M'_{n_1...n_d}$) if for $a
\in M'_{n_1...n_d}$ and any "direction" $k=1,...,d$
the intersection $span({a}^k_1,...,{a}^k_{n_k})
\cap M^{'}_{\widehat{n}_k}$
is not in general position (i.e when $a$ comes onto $M'_{n_1...n_d}$ the
topology of this intersection changes).
\end{Definition}

\begin{Theorem}
The set $M'_{n_1...n_d}$ is the discriminantal set for the system (0.3).
\end{Theorem}
As it is shown in Section 2.2 $M'_{n_1...n_d}$ is an algebraic manifold.

\begin{Definition}
{\rm The space $M_{n_{1}...n_{d}}$ is called the space of {\it inner format}
if for any  $a\in M_{n_{1}...n_{d}}$ and any "direction" $k=1,...,d$
the intersection of $GL_{n_k}$ orbit of $a$
with the submanifold $M'_{\widehat{n}_k}$ of $(d-1)-$ linear
forms of corank $1$ is not empty.\\
Otherwise the format $n_1...n_d$ is called {\it grassmanian}. The
grassmanian format $n_1...n_d$ for which $M'_{n_1...n_d}$ has
codimension 1 is called {\it boundary}.}
\end{Definition}
The term "grassmanian format" is motivated by the fact that a certain
(see Section 2.2) factorization of the space $M_{n_{1}...n_{d}}$ gives
an algebraic manifold which in case $d=2$ is the grassman
manifold $G_{n_1,n_2}$.

{\bf Example}
Consider the case of 3-linear forms of format $2\times 2\times n$ for
$n\ge 4$. These are fograssmanoundary format. For a given $a
\in M_{22n}$ take its $GL_{n}$-orbit, i.e. the set of linear
combinations $S(z):=A_1z_1+...+A_nz_n$, where
$2\times 2$ matrices $A_1,...,A_n$ are $2\times 2$ slices of $a$.
The intersection of $span(A_1,...,A_n)$ with
the submanifold $M^{'}_{22}$ of degenerate $2\times 2$ matrices
corresponds to the set of such $(z_1,...,z_n)$ that $det(S(z))=0$.
The expression $det(S(z))$ is a quadratic form on $z$. The
intersection of $span(A_1,...,A_n)$ with $M^{'}_{22}$ is described in
terms of the rank of this quadratic form. Then the notion of the corank of our
$3-$linear form $a (A_1,...,A_n)$ can be formulated in terms of
the rank of quadratic form  $det(S(z))$ as follows:
\begin{Proposition}
The corank of $2\times 2\times n$ form $a (A_1,...,A_n)$ is equal
to 1 iff the rank of quadratic form $det(S(z))$ is equal to 2.
\end{Proposition}

\subsection{Proper hyperveronese manifolds}

For a set of integers $n_1,...,n_d$ let $m_d:=n_1+...+n_d+1-d$.
Take an $n_1\times ...\times n_r \times n_{r+1}\times m_{r+1}$ matrix
$(a_{i_1..i_{r+1}j})$. Denote $T_i:=a^{(r+1)}_i,\quad i=1,...,n_{r+1}$
%Denote by $T_1,...,T_{n_{r+1}}$ matrices which are
the $n_1\times ...\times n_r\times m_{r+1}$ "slices" of $a$.
%, i.e. $$(T_i)_{i_1...i_rj}=a_{i_1...i_rij}.$$
For $(z_1,...,z_{n_{r+1}})\in {(\C^{*})^{n_{r+1}}}$ let
$S(z_1,...,z_{n_{r+1}}):=z_1T_1+...+z_{n_{r+1}}T_{n_{r+1}}$ be points of
$span(T_1,...,T_{n_{r+1}})$. Denote by $\Delta_{j_1...j_{m_r}}(S)$,
where $(j_1,...,j_{m_r})\in \left(\matrix{[m_{r+1}]\cr m_r}\right)$ the
$n_1\times ...\times n_r \times m_r$ minors of $S$.
Then the problem of finding the intersection of $span(T_1,...,T_{n_{r+1}})
\cap M'_{n_1...n_rm_{r+1}}$ is the problem of solving the system:
%(6)
\begin{equation}
\Delta_{j_1...j_{m_r}}(S)=0,
\quad (j_1,...,j_{m_r})\in \left(\matrix{[m_{r+1}]\cr m_r}\right)
\end{equation}

\begin{Proposition}
The system (6) has a nonzero solution $z\in (\C^*)^{n_{r+1}}$ iff $a$
belongs to a submanifold of codimension 1.
\end{Proposition}
According to Definitions 2.1.1 and 2.1.2 this means that the format
$n_1\times ...\times n_{r+1}\times m_{r+1}$ is boundary grassmanian.

The set $\Delta=(\Delta_{j_1...j_{m_r}})$ of minors gives us the components of
a vector in $\PP^{\left(\matrix{[m_{r+1}] \cr m_r}\right)}$. So we have
a map:
$$\phi :(C^*)^{n_{r+1}}\to \PP^{\left( \matrix{m_{r+1}\cr m_r}\right)}$$
$$(z_1,...,z_{n_{r+1}})\mapsto \Delta (z)$$
the image of which is a manifold $\cal{V}$ parametrized by
$(z_1,...,z_{n_{r+1}})$.

\begin{Theorem}
The following statements are equivalent:\\
1) $det(a_{i_1...i_{r+1}j})=0$\\
2) the manifold $\cal{V}$ is singular\\
3) the system
$$\Delta_{j_1...j_{m_r}}(S(z))=0,\quad (j_1,...,j_{m_r})\in
\left( \matrix{[m_{r+1}]\cr m_r}\right)$$
has solutions in $(\C^*)^{n_{r+1}}$ and $\phi$ gives a $1-1$ correspondence
between the solutions of the system and the singular points of $\cal{V}$.
\end{Theorem}
Compairing this statement with Teorem 1.2.1 we are lead to the following:

\begin{Definition}
{\rm For $(a)\in M_{n_1...n_{r+1}m_{r+1}}$ the manifold $\cal{V}$ is called
{\it proper hyperveronese manifold}}.
\end{Definition}

\subsection{Hypergrassmanians}
Let $(a_{i_1...i_rj})\in M_{n_1...n_rm}$ where $m>m_r=1+n_1+...+n_r-r$.
The $n_1\times ...\times n_r\times m_r$
minors $\Delta_{j_1...j_{m_r}}$ of $(a)$ are invariants of $GL_{n_1}\times
...\times GL_{n_r}$ action on $M_{n_1...n_rm}=(V_{n_1}\otimes ...\otimes
V_{n_r}\otimes V_{m})^*$.
For $1\le k\le r$ take the set $a^{(k)}_1,...,a^k_{n_k}$ of
$n_1\times ...\times {\widehat{n}_k}\times ...\times m$
"slices" of $(a)$ in $k-$th direction (see Section 2.1).
\begin{Proposition}
The intersection $span(a^{(k)}_1,...,a^{(k)}_{n_k})\cap M'_{\widehat{n}_k}$
is not empty iff
$$\Delta_{j_1...j_{m_r}}(a)=0,\quad for all (j_1,...,j_{m_r})\in
\left( \matrix {[m]\cr m_r}\right).$$
\end{Proposition}
This implies that the map
$${\cal{P}}: M_{n_1...n_rm}\to \PP^{\left( \matrix{m\cr m_r}\right)}$$
$$(a)\mapsto \Delta (a)$$
has the kernel $M'_{n_1...n_rm}$, induces an injection of the open stratum
of the space of orbits $M_{n_1...n_rm}/GL_{n_1}\times ...\times GL_{n_r}$
into $\PP^{\left( \matrix {[m]\cr m_r}\right)}$ and its image is a projective
algebraic manifold ${\cal{G}}_{n_1...n_r,m}$.
\begin{Definition}
{\rm The manifold ${\cal{G}}_{n1...n_r,m}$ is called {\it hypergrassmanian}.}
\end{Definition}
So the coordinates on a hypergrassmanian (the open stratum of the factor space
of the space $M_{n_1...n_rm}$ of forms of grassmanian format), as in the
particular case (for $r=1$) of grassman manifolds, are given by the minors
of boundary format ("hyperplukker coordinates").

%%%%%%%%%%%%%%%%%%%%%%%%%%%  SECTION 3
\section{Algorithm}
Here we give an algorithm for computing the determinants of
matrices of boundary format (``hyperplukker coordinates'').
\subsection{Basic example}
Let $P_0,...,P_{d}$ be a sequence of ordered sets
$P_k$, such that $|P_k|=k+1$.\\
{\bf Notation.} Denote by $\cal{C}$ the space of
sequences $q:=(q_1,q_2,...,q_d)$, where $q_k\in{P_k}$.

For a pair of such sequences $q^{'},q^{''}\in{\cal{C}}$ we say that
$q'\le{q^{''}}$, if $\exists K<{d}$, such that
$q^{'}_k\le{q^{''}_k}$ in $P_k$, for $k>{K}$.

\font\sixrm=cmbsy10

%%%%%  FIGURE 1

\begin{center}
\begin{picture}(400,120)(0,-120)

\put(0,-120){\makebox(100,120){
\begin{picture}(100,100)(-20,-100)
%\put(50,10){\makebox(0,0)[t]{{.}\quad{.}\quad{.}}}
\put(10,-2){\makebox(0,0)[t]{\sixrm\symbol{15}}}
\put(30,-2){\makebox(0,0)[t]{\sixrm\symbol{15}}}
\put(50,-2){\makebox(0,0)[t]{\sixrm\symbol{15}}}
\put(70,-2){\makebox(0,0)[t]{\sixrm\symbol{15}}}
\put(90,-2){\makebox(0,0)[t]{\sixrm\symbol{15}}}
\put(20,-22){\makebox(0,0)[t]{\sixrm\symbol{15}}}
\put(40,-22){\makebox(0,0)[t]{\sixrm\symbol{15}}}
\put(60,-22){\makebox(0,0)[t]{\sixrm\symbol{15}}}
\put(80,-22){\makebox(0,0)[t]{\sixrm\symbol{15}}}
\put(30,-42){\makebox(0,0)[t]{\sixrm\symbol{15}}}
\put(50,-42){\makebox(0,0)[t]{\sixrm\symbol{15}}}
\put(70,-42){\makebox(0,0)[t]{\sixrm\symbol{15}}}
\put(40,-62){\makebox(0,0)[t]{\sixrm\symbol{15}}}
%\put(40,-62){\makebox(0,0)[t]{\sevenrm\symbol{35}}}
%\put(40,-62){\makebox(0,0)[t]{\times}}
\put(60,-62){\makebox(0,0)[t]{\sixrm\symbol{15}}}
\put(50,-82){\makebox(0,0)[t]{\sixrm\symbol{15}}}
\put(-20,0){\makebox(0,0)[t]{$P_4$}}
\put(-10,-20){\makebox(0,0)[t]{$P_3$}}
\put(0,-40){\makebox(0,0)[t]{$P_2$}}
\put(10,-60){\makebox(0,0)[t]{$P_1$}}
\put(20,-80){\makebox(0,0)[t]{$P_0$}}

\put(40,-110){\makebox(20,10){Fig.1}}
\end{picture}}}

%%%%%  FIGURE 2

\put(150,-120){\makebox(100,120){
\begin{picture}(100,100)(-20,-100)
%\put(50,10){\makebox(0,0)[t]{{.}\quad{.}\quad{.}}}
\put(10,-2){\makebox(0,0)[t]{\sixrm\symbol{15}}}
\put(30,-2){\makebox(0,0)[t]{\sixrm\symbol{2}}}
\put(50,-2){\makebox(0,0)[t]{\sixrm\symbol{15}}}
\put(70,-2){\makebox(0,0)[t]{\sixrm\symbol{15}}}
\put(90,-2){\makebox(0,0)[t]{\sixrm\symbol{15}}}
\put(20,-22){\makebox(0,0)[t]{\sixrm\symbol{15}}}
\put(40,-22){\makebox(0,0)[t]{\sixrm\symbol{2}}}
\put(60,-22){\makebox(0,0)[t]{\sixrm\symbol{15}}}
\put(80,-22){\makebox(0,0)[t]{\sixrm\symbol{15}}}
\put(30,-42){\makebox(0,0)[t]{\sixrm\symbol{15}}}
\put(50,-42){\makebox(0,0)[t]{\sixrm\symbol{15}}}
\put(70,-42){\makebox(0,0)[t]{\sixrm\symbol{2}}}
\put(40,-62){\makebox(0,0)[t]{\sixrm\symbol{2}}}
%\put(40,-62){\makebox(0,0)[t]{\sevenrm\symbol{35}}}
%\put(40,-62){\makebox(0,0)[t]{\times}}
\put(60,-62){\makebox(0,0)[t]{\sixrm\symbol{15}}}
\put(50,-82){\makebox(0,0)[t]{\sixrm\symbol{15}}}

\put(40,-110){\makebox(20,10){Fig.2}}
\end{picture}}}

%%%%%  FIGURE 3

\put(300,-120){\makebox(100,120){
\begin{picture}(100,100)(-20,-100)
%\put(50,10){\makebox(0,0)[t]{{.}\quad{.}\quad{.}}}
\put(10,-2){\makebox(0,0)[t]{\sixrm\symbol{15}}}
\put(30,-2){\makebox(0,0)[t]{\sixrm\symbol{2}}}
\put(50,-2){\makebox(0,0)[t]{\sixrm\symbol{15}}}
\put(70,-2){\makebox(0,0)[t]{\sixrm\symbol{15}}}
\put(90,-2){\makebox(0,0)[t]{\sixrm\symbol{15}}}
\put(20,-22){\makebox(0,0)[t]{\sixrm\symbol{15}}}
\put(40,-22){\makebox(0,0)[t]{\sixrm\symbol{2}}}
\put(60,-22){\makebox(0,0)[t]{\sixrm\symbol{15}}}
\put(80,-22){\makebox(0,0)[t]{\sixrm\symbol{15}}}
\put(30,-42){\makebox(0,0)[t]{\sixrm\symbol{15}}}
\put(50,-42){\makebox(0,0)[t]{\sixrm\symbol{15}}}
\put(70,-42){\makebox(0,0)[t]{\sixrm\symbol{2}}}
\put(40,-62){\makebox(0,0)[t]{\sixrm\symbol{2}}}
%\put(40,-62){\makebox(0,0)[t]{\sevenrm\symbol{35}}}
%\put(40,-62){\makebox(0,0)[t]{\times}}
\put(60,-62){\makebox(0,0)[t]{\sixrm\symbol{15}}}
\put(50,-82){\makebox(0,0)[t]{\sixrm\symbol{15}}}

\thicklines\put(51,-82){\vector(1,2){7}}
\thicklines\put(39,-62){\vector(-1,2){7}}
\thicklines\put(59,-62){\vector(-1,2){7}}
\thicklines\put(29,-42){\vector(-1,2){7}}
\thicklines\put(51,-42){\vector(1,2){7}}
\thicklines\put(71,-42){\vector(1,2){7}}
\thicklines\put(19,-22){\vector(-1,2){7}}
\thicklines\put(41,-22){\vector(1,2){7}}
\thicklines\put(61,-22){\vector(1,2){7}}
\thicklines\put(81,-22){\vector(1,2){7}}

\put(40,-110){\makebox(20,10){Fig.3}}
\end{picture}}}
%\put(2.5,-4){\mbox(0,0){\nearrow}}
\end{picture}
\end{center}

The sequence $P_0,...,P_4$ can be represented
in the form of a diagram on Fig.1, where $(k+1)-$th row represents
the elements of the set $P_k$ ordered from the left to the right.
A sequence
$q$ may be represented in the form of a diagram on Fig.2, where
$q_k$ is represented by a cross in the $(k+1)-$th row.

For a given $q_{k+1}\in{P_{k+1}}$ there is a unique
ordering-preserving injection
$f^q_{k}:P_k\hookrightarrow{P_{k+1}}$, such that ${q_{k+1}}${\put(5,3)
{\makebox(0,0){$\in$}}}{\put(0,3){\makebox(0,0){$\not$}}}
\quad${f^q_{k}(P_k)}$ and vise versa.

For $q$ represented by Fig.2 the corresponding sequence
$f^q:={(f^q_1,...,f^q_4)}$ of injections may be represented in
form of a diagram on Fig.3.

Then each sequence $q$ defines a sequence $p(q):=(p_1,...,p_d)$, where
$p_k\in{P_k}$, such that $p_{k+1}=f^q_{k}(p_k)$. The sequence
$p$, corresponding to the sequence $q$ on Fig.2, may be represented
in form of a diagram on Fig.4, where the path goes through the
elements $p_k$. The rule for drawing a path, corresponding to a
given sequence $q$ may be formulated as follows:

\begin{em}
if $q_{k+1}\in{P_{k+1}}$ lies (on the diagram) to the left from $p_k$,
then we have to go from $p_k$ to the right closest to it point
of $P_{k+1}$, if $q_{k+1}\in{P_{k+1}}$ lies to the right from
$p_k$, then we have to go from $p_k$ to the left closest
point of $P_{k+1}$.
\end{em}

On the other hand, for a given $q_k$ the corresponding
$f^q_k:P_k\to{P_{k+1}}$ defines a partition of $P_k$ into two
parts $P^{-}_{k}(q)$ and $P^{+}_{k}(q)$, where
$$P^{-}_{k}(q):=\{{p_k\in{P_k}|f^q_k(p_k)<q_{k+1}}\}$$
$$P^{+}_{k}(q):=\{{p_k\in{P_k}|f^q_k(p_k)>q_{k+1}}\}$$
The sequence of partitions, corresponding to the sequence of
ordering pre- serving injections
$f^q=(f^q_0,f^q_1,f^q_2,f^q_3)$ on Fig.3, may be represented in form of
a diagram on Fig.5, where in the $k-$th row the elements
of $P^{-}_{k}(q)$ are represented by characters ``1'' and the
the elements of $P^{+}_{k}(q)$ are represented by characters ``2''.

%%%%%% FIGURE 4

\begin{center}
\begin{picture}(400,120)(0,-120)
\put(-20,-120){\makebox(100,120){
\begin{picture}(100,100)(-20,-100)
\put(10,-2){\makebox(0,0)[t]{\sixrm\symbol{15}}}
\put(30,-2){\makebox(0,0)[t]{\sixrm\symbol{2}}}
\put(50,-2){\makebox(0,0)[t]{\sixrm\symbol{15}}}
\put(70,-2){\makebox(0,0)[t]{\sixrm\symbol{15}}}
\put(90,-2){\makebox(0,0)[t]{\sixrm\symbol{15}}}
\put(20,-22){\makebox(0,0)[t]{\sixrm\symbol{15}}}
\put(40,-22){\makebox(0,0)[t]{\sixrm\symbol{2}}}
\put(60,-22){\makebox(0,0)[t]{\sixrm\symbol{15}}}
\put(80,-22){\makebox(0,0)[t]{\sixrm\symbol{15}}}
\put(30,-42){\makebox(0,0)[t]{\sixrm\symbol{15}}}
\put(50,-42){\makebox(0,0)[t]{\sixrm\symbol{15}}}
\put(70,-42){\makebox(0,0)[t]{\sixrm\symbol{2}}}
\put(40,-62){\makebox(0,0)[t]{\sixrm\symbol{2}}}
%\put(40,-62){\makebox(0,0)[t]{\sevenrm\symbol{35}}}
%\put(40,-62){\makebox(0,0)[t]{\times}}
\put(60,-62){\makebox(0,0)[t]{\sixrm\symbol{15}}}
\put(50,-82){\makebox(0,0)[t]{\sixrm\symbol{15}}}

\thicklines\put(51,-82){\vector(1,2){7}}
\thicklines\put(59,-62){\vector(-1,2){7}}
\thicklines\put(51,-42){\vector(1,2){7}}
\thicklines\put(61,-22){\vector(1,2){7}}
\put(40,-110){\makebox(20,10){Fig.4}}
\end{picture}}}
%\end{picture}
%\end{center}

%%%%  FIGURE 5

\put(130,-120){\makebox(100,120){
\begin{picture}(100,100)(-20,-100)
%\put(50,10){\makebox(0,0)[t]{{.}\quad{.}\quad{.}}}
\put(10,-2){\makebox(0,0)[t]{\sixrm\symbol{15}}}
\put(30,-2){\makebox(0,0)[t]{\sixrm\symbol{2}}}
\put(50,-2){\makebox(0,0)[t]{\sixrm\symbol{15}}}
\put(70,-2){\makebox(0,0)[t]{\sixrm\symbol{15}}}
\put(90,-2){\makebox(0,0)[t]{\sixrm\symbol{15}}}
\put(20,-22){\makebox(0,0)[t]{1}}
\put(40,-22){\makebox(0,0)[t]{2}}
\put(60,-22){\makebox(0,0)[t]{2}}
\put(80,-22){\makebox(0,0)[t]{2}}
\put(30,-42){\makebox(0,0)[t]{1}}
\put(50,-42){\makebox(0,0)[t]{2}}
\put(70,-42){\makebox(0,0)[t]{2}}
\put(40,-62){\makebox(0,0)[t]{1}}
\put(60,-62){\makebox(0,0)[t]{1}}
\put(50,-82){\makebox(0,0)[t]{2}}
%\put(-20,0){\makebox(0,0)[t]{$P_5$}}
%\put(-10,-20){\makebox(0,0)[t]{$P_4$}}
%\put(0,-40){\makebox(0,0)[t]{$P_3$}}
%\put(10,-60){\makebox(0,0)[t]{$P_2$}}
%\put(20,-80){\makebox(0,0)[t]{$P_1$}}

\put(40,-110){\makebox(20,10){Fig.5}}
\end{picture}}}

%%%%   FIGURE 6

\put(280,-120){\makebox(100,120){
\begin{picture}(100,100)(-20,-100)
%\put(50,10){\makebox(0,0)[t]{{.}\quad{.}\quad{.}}}
\put(10,-2){\makebox(0,0)[t]{1}}
\put(30,-2){\makebox(0,0)[t]{2}}
\put(50,-2){\makebox(0,0)[t]{3}}
\put(70,-2){\makebox(0,0)[t]{4}}
\put(90,-2){\makebox(0,0)[t]{5}}
\put(20,-22){\makebox(0,0)[t]{1}}
\put(40,-22){\makebox(0,0)[t]{2}}
\put(60,-22){\makebox(0,0)[t]{2}}
\put(80,-22){\makebox(0,0)[t]{2}}
\put(30,-42){\makebox(0,0)[t]{1}}
\put(50,-42){\makebox(0,0)[t]{2}}
\put(70,-42){\makebox(0,0)[t]{2}}
\put(40,-62){\makebox(0,0)[t]{1}}
\put(60,-62){\makebox(0,0)[t]{1}}
\put(50,-82){\makebox(0,0)[t]{2}}
%\put(-20,0){\makebox(0,0)[t]{$P_5$}}
%\put(-10,-20){\makebox(0,0)[t]{$P_4$}}
%\put(0,-40){\makebox(0,0)[t]{$P_3$}}
%\put(10,-60){\makebox(0,0)[t]{$P_2$}}
%\put(20,-80){\makebox(0,0)[t]{$P_1$}}

\thicklines\put(52,-81){\vector(1,2){5}}
\thicklines\put(58,-61){\vector(-1,2){5}}
\thicklines\put(52,-41){\vector(1,2){5}}
\thicklines\put(62,-21){\vector(1,2){5}}

\put(40,-110){\makebox(20,10){Fig.6}}
\end{picture}}}

\end{picture}
\end{center}
%\font\sixrm=sixrm

{\bf Notation.}Denote by $j(q)$ the ordinal number of the element
$p_{d}(q)\in{P_{d}}$ with respect to the ordering on $P_{d}$.
For a given pair $q^{'},q^{''}\in{\cal{C}}$ denote by
$[i_1,...,i_d,j](q^{'},q^{''})$ the sequence of indecies
$(i_1,...,i_d,j)$ such that $$j=j(Q)
%: $$i(q^{'},q^{''}):=[i_1,{...},i_d](q^{'},q^{''}){\quad}
{\quad}and{\quad}i_{k+1}={\left\{\begin{array}{ll}
1,\quad{{\rm if \ }p_{k}(q^{'})\in{P^{-}_{k}(q^{''})}}\\
2,\quad{{\rm if \ }p_{k}(q^{'})\in{P^{+}_{k}(q^{''})}}
\end{array}\right.}$$
For example, for $q^{'}=q^{''}=q$ from Fig.2 we can write $(i_1,...,i_4,j)$
using the diagram on Fig.6, which may be viewed as a ``superposition''
of Fig.4 and Fig.5, $(i_1,...,i_4,j)=[i_1,i_2,i_3,i_4,j](q,q)=(2,1,2,2,4)$.

\begin{Proposition}\label{P1}
$$\prod^{}_{q\in{\cal{C}}}{a_{[i_1,...,i_d,j](q,q)}}$$
is a monomial of discriminant of polylinear form of
dimension ${\underbrace{2\times{2}\times{...}\times{2}}_{\mbox
{\scriptsize{d\ {times}}}}}$
$\times{(d+1)}$ with coefficients $(a_{i_1,...,i_d,j})$.
\end{Proposition}

This monomial will be called {\it diagonal{\ }monomial\/}.

%%%%   SECTION 2
\subsection{General algorithm}
%~
Let $n_1,...,n_d$ be a sequense of positive integers and set $n_0=1$.
Denote $m_k:=n_0+n_1+...+n_k-k$.
Let $P_0,...,P_d$ be a sequence of ordered sets, such that
$|P_k|=m_k$. Let $T_1,...,T_d$ be a sequence of ordered sets,
such that $|T_k|=n_k$.\\
{\bf Notation.} Denote by ${\cal{C}}_k:=\{Q_k\subset P_k\
:\ |Q_k|=m_k-m_{k-1}\}$ and by $\cal{C}$ the space of sequences
$Q:=(Q_1,...,Q_d)$, where $Q_k\in {\cal{C}}_k$.
For a pair
$S^{'},S^{''}\subset{P_k}$, where $S^{'}=(s^{'}_1,...,s^{'}_{n_k}),\
S^{''}=(s^{''}_1,...,s^{''}_{n_k})$ and the elements of $S^{'}$ and
$S^{''}$ are written in the order induced by the ordering on $P_k$,
we say that $S^{'}\le{S^{''}}$ if $\exists N<n_k$, such that $s^{'}
_n\le{s^{''}_n}$ for $n>N$. This gives the ordering of ${\cal{C}}_k$.
For a pair of sequences $Q^{'},Q^{''}\in{\cal{C}}$ we say that
$Q^{'}\le{Q^{''}}$ if $\exists K<d$, such that $Q^{'}_k\le{Q^{''}_k}$,
for $k>K$.
For a given $Q_{k+1}\in {\cal{C}}_{k+1}$ there exists a unique
order preserving injection $f^Q_k:P_k\hookrightarrow P_{k+1}\setminus
Q_{k+1}$
Then each sequence
$Q\in{\cal{C}}$ defines a sequence $p(Q):=(p_0,p_1,...,p_d)$, where
$p_k\in P_k$, such that $p_{k+1}=f^Q_k(p_k)$ which is caled {\em $Q-$path}.
An example of a $Q-$path is shown on Fig.4.

On the other hand, since $P_k$ are ordered, then a given $Q_k$
defines a partition $R_k=(R^1_k,...,R^{n_k}_k)$ of $P_k\setminus Q_k$
into $n_k$ subsets $R^i_k$.
%~
Then any injection $\phi : P_{k-1}\hookrightarrow P_{k}\setminus
Q_k$ induce a map $g_k: P_{k-1}\to T_k$ as follows:
%(7)
\begin{equation}
if p\in \phi^{-1}(R^i_k)\quad then g_k(p)=i\in T_k
\end{equation}
where we write number $i$ instead of the $i-$th element of $T_k$.
\begin{Definition}
{\rm For a given $Q_k\in {\cal{C}}_k$ such a map will be called
{\it $Q-$admissible} and the sequence $g^Q=(g_1,...,g_d)$ of $Q-$admissible
maps is called a {\it $Q-$diagram}. If in a $Q-$diagram all $g_k$ are induced
by
the order preserving injections $f^Q_k$, then this $Q-$diagram is called
{\it initial}}.
\end{Definition}
An example of the initial $Q-$diagram is shown on Fig.5.
Now let us take a pair $Q^{'},Q^{''}\in\cal{C}$. Then for $Q^{'}$
we take the $Q^{'}-$path $p(Q^{'})=(p_0,p_1,...,p_d)$ and for $Q^{''}$
we choose a $Q^{''}-$diagram $g^{Q^{''}}=(g_1,...,g_d)$ from the set of
$Q^{''}-$admissible diagrams. This gives us a sequence of indecies:
%(8)
\begin{equation}
I(Q^{'},g^{Q^{''}})=(i_1,...,i_d), \quad where
i_k:=g^{Q^{''}}(p^{Q^{'}}_{k-1}).
\end{equation}
\begin{Definition}
{\rm The pair $(p^{Q^{'}},g^{Q^{''}})$ is called {\it $Q^{'}-$path over a
$Q^{''}-$diagram}.}
\end{Definition}
An example of this construction is shown on Fig.6.
{\bf Notation.} Denote by ${\cal{C}}^{(k)}$ the space of subsequences
$Q^{(k)}:=(Q_{k+1},...,Q_{d})$.
\begin{Definition}
{\rm For a given $Q^{(k)}\in{{\cal{C}}^{(k)}}$ two sequences
$Q^{'},Q^{''}\in{\cal{C}}$ will be called {\em $Q^{(k)}-$conjugate}
if $Q^{'}_m=Q^{''}_m=Q^{(k)}_m$ for $m>k$.
For a given $Q^{(k)}\in{{\cal{C}}^{(k)}}$ and $p:=(p_1,...,p_d)$
two sequences $Q^{'},Q^{''}
\in{\cal{C}}$ will be called {\it $(Q^{(k)},p_k)-$conjugate}, if
$Q^{'}_m=Q^{''}_m=Q^{(k)}_m$ for $m>k$, and $p(Q^{'})_k=
p(Q^{''})_k=p_k\in{P_k}$.}
\end{Definition}
{\bf Notation.} Denote by ${\cal{D}}{(Q^{(k)})}$ the set of all
$Q^{(k)}-$conjugate sequences and by ${\cal{D}}{(Q^{(k)},p_k)}$
the set of all $(Q^{(k)},p_k)-$conjugate sequences in $\cal{C}$.
Set by definition ${\cal{D}}(Q^{(d)}):=\cal{C}$.
%
%     PROPOSITION
%
\begin{Proposition}
For a given $k\le d$ $$|{\cal{D}}{(Q^{(k)},p_k)}|=\frac{(m_k-1)!}
{(n_1-1)!...(n_k-1)!}$$  for any $p_k\in P_k$.
\end{Proposition}

On each ${\cal{D}}{(Q^{(k)},p_k)}$ there is an ordering induced by
the ordering on $\cal{C}$.\\
{\bf Notation.} For a given $k$ denote by ${\cal{L}}_{k}$
the set of ordinal numbers,
%such that $|{\cal{L}}_{k}|=
ennumerating the elements of each ${\cal{D}}{(Q^{(k)},p_k)}$.

{\em Remark:} Proposition 2.1 implies, that
$|{\cal{L}}_{0}|=1$ and $|{\cal{L}}_{1}|=1$.\\
Then for a given $Q\in {\cal{C}}$ we
have a sequence of integer numbers $L(Q):=(l_1,...,
l_d)$, where $l_k (Q)$ is the ordinal number of $Q$
in ${\cal{D}}(Q^{(k)},p_k)$.\\
{\bf Notation.} For a given $Q^{(k)}$ and $l\in
{\cal{L}}_{k-1}$ denote by ${\cal{E}}(Q^{(k)},l)$ the
set of sequences $Q$, such that $l_{k-1}(Q)=l$.
Denote by $S_{{\cal{C}}_k}$ the group of permutations of the elements of
${\cal{C}}_k$.

Then $S_{{\cal{C}}_k}$ acts on ${\cal{E}}(Q^{(k)},l)$ as follows:
for $\sigma\in S_{{\cal{C}}_k}$ and a given $Q=(Q_1,...,Q_k,Q^{(k)})
\in {\cal{E}}(Q^{(k)},l)$, $$\sigma Q=(Q_1,...,
\sigma Q_k,Q^{(k)}).$$
If to each ${\cal{E}}(Q^{(k)},l_{k-1})$ we put into correspondence a group
$S(Q^{(k)},l_{k-1})\cong S_{{\cal{C}}_k}$ with the action discribed above,
then on the whole ${\cal{C}}$ we have the action of the group
$$\Sigma:=\prod^{d-1}_{k=0}\ \prod^{}_{Q^{(k+1)}}
\ \prod^{}_{l_{k}}{S(Q^{(k+1)},l_{k})}$$
Each $Q\in{\cal{C}}$ defines the following subgroup $\Sigma_Q$ of
$\Sigma$:$$\Sigma_Q=S_{{\cal{C}}_1}(Q)\times
...\times{S_{{\cal{C}}_d}}(Q)$$ where $S_{{\cal{C}}_k}(Q)=S(Q^{(k)},
\l_{k-1}(Q))$. $\Sigma_Q$ acts on $Q$ componentwise: for $\tau=
(\tau_1,...,\tau_d)\in \Sigma_Q$
$$\tau Q=(\tau_{1}Q_1,...,\tau_d Q_d).$$
The role of the group $\Sigma$ is analogous to the role of
symmetric group in calculation of the determinant of $n\times n$ matrix.

%%%           PROPOSITION
\begin{Proposition}
If $Q\in\cal{C}$, $\sigma\in\Sigma$ and $\sigma_Q\in\Sigma_Q\subset
\Sigma$ is the $\Sigma_Q$-component of $\sigma$ then
$$\sigma Q=\sigma_Q Q$$.
\end{Proposition}

For $\sigma\in \Sigma$ denote $$sign(\sigma):=\prod^{d-1}_{k=0}
\ \prod^{}_{Q^{(k+1)}}
\ \prod^{}_{l_{k}}{sign(\sigma(Q^{(k+1)},l_{k}))}$$ where
$sign(\sigma(Q^{(k+1)},l_{k}))$ is the signature of permutation
$\sigma(Q^{(k+1)},l_{k})\in S(Q^{(k+1)},l_{k})$.

Let us put into correspondence to each $Q^{(k)}\in {\cal{C}}^{(k)}$
and $l_k\in{\cal{L}}_k$ a set $G(Q^{(k)},l_k)$
of $Q_{k+1}$-admissible functions $g_{k+1}$ and take their direct product:
$$\Gamma:=\prod^{d-1}_{k=0}\prod^{}_{Q^{(k)}}\prod^{}_{l_{k}}{G(Q^{(k)},l_k)}$$
Denote by $j(Q)$ the ordinal number of the element $p_{d}(Q)\in P_d$
with respect to the ordering on $P_d$.
For given $Q\in\cal{C}$, $\sigma\in\Sigma$ and $\gamma\in\Gamma$ denote
by $[i_1,...,i_d](Q,\sigma Q,\gamma)$ the sequence of indecies
$(i_1,...,i_d,j)$, where $$j=j(Q),{\quad} i_k=g_k(p_{k-1}(Q)){\quad}
and{\quad}g_k=\gamma(Q^{(k)},l_{k}(Q))\in{G((\sigma Q)^{(k)},l_k(Q))}.$$
%
%   THEOREM
%
\begin{Theorem}
If $\Omega=(a_{i_1...,i_dj})$ is a polylinear form of
dimension $n_1\times{...}\times{n_d}\times{m_d}$, then its discriminant
$$D_{\Omega}=\sum^{}_{\sigma\in\Sigma}{sign(\sigma)}{\sum^{}_{\gamma
\in\Gamma}}\ {\prod^{}_{Q\in\cal{C}}}{a_{[i_1,...,i_d,j]
(Q,\sigma Q,\gamma)}}$$
\end{Theorem}
{\bf Example.} {\em Let $d=1$, $n_1=n$}. \\Then $m_1=n$,
${\cal{C}}={\cal{C}}_1=\{{Q_1\subset P_1:\ |Q|=n-1}\}$ and each $Q$ is
defined by the value of $j(Q)$. Since
$|{\cal{L}}_0|=1$, then for each $Q=Q_1$ the set $G(Q^{(0)})$
consists of only one element $g^Q$, such that the ordinal number
of $g^{Q}(p_0)$ in $T_1$ is equal to $j(Q)$, and $\Gamma=\prod^{}
_{Q^{(0)}}G(Q^{(0)})$ consists of only one element $\gamma=\prod^{}_
{Q\in\cal{C}}g^{Q}$. Also $\Sigma=S_{{\cal{C}}_1}\cong S_n$.
Then for $\Omega=(a_{ij})_{1\le{i,j}\le{n}}$
$$D_{\Omega}=\sum^{}_{\sigma\in S_n}{sign(\sigma)}\prod^{}_{Q\in\cal{C}}
{a_{[i,j](Q,\sigma Q,\gamma)}}=
\sum^{}_{\sigma\in S_n}{sign(\sigma)}\prod^{n}_{j=1}{a_{\sigma
j,j}}$$ is the determinant of the square matrix $(a_{ij})_
{1\le{i,j}\le{n}}$.

%%%%%%%%%%%%%%%% SECTION 4
\section{Closed determinant}
\begin{Definition}
{\rm For $n_1\times {...}\times n_d$ $d-$dimensional
matrix $(a_{i_1...i_d})$ the product of all its minors
(including the determinant) is called
the {\it closed determinant} (denoted by $Det(a)$)}
\end{Definition}
The term "closed" comes from the fact that as an algebraic manifold $Det(a)$
corresponds to the projectively dual to the closure of the
$(\C^*)^{n_1+...+n_d}$ orbit of $(1,...,1)\otimes {...}\otimes (1,...,1)\in
V_{n_1}\otimes{...}\otimes V_{n_d}$.

{\bf Example.} {\em $2\times 2\times 2$ matrix}.
Let $(a_{i_1i_2i_3})_{i_1,i_2,i_3=1,2}$ be a $2\times 2\times 2$ matrix.
Then its closed determinant
%(9)
%^%\begin{equation}
$$Det(a)=a_{111}a_{112}a_{121}a_{122}a_{211}a_{212}a_{221}a_{222}\times$$
$$\times (a_{111}a_{122}-a_{121}a_{112})(a_{211}a_{222}-a_{221}a_{212})
(a_{111}a_{212}-a_{211}a_{112})(a_{121}a_{222}-a_{221}a_{122})
(a_{111}a_{221}-a_{211}a_{121})(a_{112}a_{222}-a_{212}a_{122})\times$$
$$\times
(a^2_{111}a^2_{222}+a^2_{112}a^2_{221}+a^2_{121}a^2_{212}+a_{211}^2a_{122}^2
-2a_{111}a_{121}a_{212}a_{222}-2a_{111}a_{211}a_{122}a_{222}
-2a_{111}a_{112}a_{221}a_{222}-2a_{121}a_{221}a_{112}a_{212}
-2a_{211}a_{221}a_{112}a_{122}-2a_{212}a_{211}a_{121}a_{122}
+4a_{111}a_{221}a_{212}a_{122}+4a_{121}a_{211}a_{112}a_{222})$$
%^%\end{equation}
\begin{Proposition}
The degree of the closed determinant of $d-$dimensional matrix of format
$n_1\times {...}\times n_d$ is equal to the degree of the (ordinary)
determinant of $(d+1)-$dimensional matrix of boundary format $n_1\times
{...}\times n_d\times (1+n_1+...+n_d-d)$.
\end{Proposition}
Let us take for each initial $Q-$diagram the $Q-$path on it (see Definition
3.2.1).
The corresponding set of indecies will be denoted by $I(Q)$.
\begin{Theorem}
Let $(a_{i_1...i_d})$ be a $d-$dimensional matrix.
The monomial
\begin{equation}
\prod^{}_{Q\in \cal C}a_{I(Q)}
\end{equation}
is a monomial of the closed determinant $Det(a)$.
\end{Theorem}
\begin{Definition}
{\rm The monomial} $\prod^{}_{Q\in {\cal C}}a_{I(Q)}$ {\rm is called {\it
diagonal}}.
\end{Definition}
There is an algorithm given in terms of paths over $Q-$diagrams (of which
the evidencies of existence are Proposition 4.1 and Theorem 4.1) of
computing closed determinants of polylnear forms, which will be published
separately. Here we give an example of this procedure.\\
{\bf Example.} {\em $2\times 2\times 2$ matrix}.
As a tool for our calculation let us draw the set of initial $Q-$diagrams
with corresponding $Q-$paths on them

%%%  FIGURE 7

\begin{center}
\begin{picture}(380,300)(0,-300)

%% level 1

 \put(0,-50){\makebox(40,50){
  \begin{picture}(40,50)(0,-50)
\put(0,-2){\makebox(0,0)[t]{1}}
\put(20,-2){\makebox(0,0)[t]{1}}
\put(40,-2){\makebox(0,0)[t]{1}}
\put(10,-22){\makebox(0,0)[t]{1}}
\put(30,-22){\makebox(0,0)[t]{1}}
\put(20,-42){\makebox(0,0)[t]{1}}
\thicklines\put(18,-41){\vector(-1,2){5}}
\thicklines\put(8,-21){\vector(-1,2){5}}
  \end{picture}}}

\put(60,-50){\makebox(40,50){
  \begin{picture}(40,50)(0,-50)
\put(0,-2){\makebox(0,0)[t]{1}}
\put(20,-2){\makebox(0,0)[t]{1}}
\put(40,-2){\makebox(0,0)[t]{1}}
\put(10,-22){\makebox(0,0)[t]{1}}
\put(30,-22){\makebox(0,0)[t]{1}}
\put(20,-42){\makebox(0,0)[t]{2}}
\thicklines\put(22,-41){\vector(1,2){5}}
\thicklines\put(28,-21){\vector(-1,2){5}}
  \end{picture}}}

 \put(140,-50){\makebox(40,50){
  \begin{picture}(40,50)(0,-50)
\put(0,-2){\makebox(0,0)[t]{1}}
\put(20,-2){\makebox(0,0)[t]{1}}
\put(40,-2){\makebox(0,0)[t]{1}}
\put(10,-22){\makebox(0,0)[t]{1}}
\put(30,-22){\makebox(0,0)[t]{2}}
\put(20,-42){\makebox(0,0)[t]{1}}
\thicklines\put(18,-41){\vector(-1,2){5}}
\thicklines\put(8,-21){\vector(-1,2){5}}
  \end{picture}}}

\put(200,-50){\makebox(40,50){
  \begin{picture}(40,50)(0,-50)
\put(0,-2){\makebox(0,0)[t]{1}}
\put(20,-2){\makebox(0,0)[t]{1}}
\put(40,-2){\makebox(0,0)[t]{1}}
\put(10,-22){\makebox(0,0)[t]{1}}
\put(30,-22){\makebox(0,0)[t]{2}}
\put(20,-42){\makebox(0,0)[t]{2}}
\thicklines\put(22,-41){\vector(1,2){5}}
\thicklines\put(32,-21){\vector(1,2){5}}
  \end{picture}}}

 \put(280,-50){\makebox(40,50){
  \begin{picture}(40,50)(0,-50)
\put(0,-2){\makebox(0,0)[t]{1}}
\put(20,-2){\makebox(0,0)[t]{1}}
\put(40,-2){\makebox(0,0)[t]{1}}
\put(10,-22){\makebox(0,0)[t]{2}}
\put(30,-22){\makebox(0,0)[t]{2}}
\put(20,-42){\makebox(0,0)[t]{1}}
\thicklines\put(18,-41){\vector(-1,2){5}}
\thicklines\put(12,-21){\vector(1,2){5}}
  \end{picture}}}

\put(340,-50){\makebox(40,50){
  \begin{picture}(40,50)(0,-50)
\put(0,-2){\makebox(0,0)[t]{1}}
\put(20,-2){\makebox(0,0)[t]{1}}
\put(40,-2){\makebox(0,0)[t]{1}}
\put(10,-22){\makebox(0,0)[t]{2}}
\put(30,-22){\makebox(0,0)[t]{2}}
\put(20,-42){\makebox(0,0)[t]{2}}
\thicklines\put(22,-41){\vector(1,2){5}}
\thicklines\put(32,-21){\vector(1,2){5}}
  \end{picture}}}

%% level 2

 \put(0,-130){\makebox(40,50){
  \begin{picture}(40,50)(0,-50)
\put(0,-2){\makebox(0,0)[t]{1}}
\put(20,-2){\makebox(0,0)[t]{1}}
\put(40,-2){\makebox(0,0)[t]{2}}
\put(10,-22){\makebox(0,0)[t]{1}}
\put(30,-22){\makebox(0,0)[t]{1}}
\put(20,-42){\makebox(0,0)[t]{1}}
\thicklines\put(18,-41){\vector(-1,2){5}}
\thicklines\put(8,-21){\vector(-1,2){5}}
  \end{picture}}}

\put(60,-130){\makebox(40,50){
  \begin{picture}(40,50)(0,-50)
\put(0,-2){\makebox(0,0)[t]{1}}
\put(20,-2){\makebox(0,0)[t]{1}}
\put(40,-2){\makebox(0,0)[t]{2}}
\put(10,-22){\makebox(0,0)[t]{1}}
\put(30,-22){\makebox(0,0)[t]{1}}
\put(20,-42){\makebox(0,0)[t]{2}}
\thicklines\put(22,-41){\vector(1,2){5}}
\thicklines\put(28,-21){\vector(-1,2){5}}
  \end{picture}}}

 \put(140,-130){\makebox(40,50){
  \begin{picture}(40,50)(0,-50)
\put(0,-2){\makebox(0,0)[t]{1}}
\put(20,-2){\makebox(0,0)[t]{1}}
\put(40,-2){\makebox(0,0)[t]{2}}
\put(10,-22){\makebox(0,0)[t]{1}}
\put(30,-22){\makebox(0,0)[t]{2}}
\put(20,-42){\makebox(0,0)[t]{1}}
\thicklines\put(18,-41){\vector(-1,2){5}}
\thicklines\put(8,-21){\vector(-1,2){5}}
  \end{picture}}}

\put(200,-130){\makebox(40,50){
  \begin{picture}(40,50)(0,-50)
\put(0,-2){\makebox(0,0)[t]{1}}
\put(20,-2){\makebox(0,0)[t]{1}}
\put(40,-2){\makebox(0,0)[t]{2}}
\put(10,-22){\makebox(0,0)[t]{1}}
\put(30,-22){\makebox(0,0)[t]{2}}
\put(20,-42){\makebox(0,0)[t]{2}}
\thicklines\put(22,-41){\vector(1,2){5}}
\thicklines\put(32,-21){\vector(1,2){5}}
  \end{picture}}}

 \put(280,-130){\makebox(40,50){
  \begin{picture}(40,50)(0,-50)
\put(0,-2){\makebox(0,0)[t]{1}}
\put(20,-2){\makebox(0,0)[t]{1}}
\put(40,-2){\makebox(0,0)[t]{2}}
\put(10,-22){\makebox(0,0)[t]{2}}
\put(30,-22){\makebox(0,0)[t]{2}}
\put(20,-42){\makebox(0,0)[t]{1}}
\thicklines\put(18,-41){\vector(-1,2){5}}
\thicklines\put(12,-21){\vector(1,2){5}}
  \end{picture}}}

\put(340,-130){\makebox(40,50){
  \begin{picture}(40,50)(0,-50)
\put(0,-2){\makebox(0,0)[t]{1}}
\put(20,-2){\makebox(0,0)[t]{1}}
\put(40,-2){\makebox(0,0)[t]{2}}
\put(10,-22){\makebox(0,0)[t]{2}}
\put(30,-22){\makebox(0,0)[t]{2}}
\put(20,-42){\makebox(0,0)[t]{2}}
\thicklines\put(22,-41){\vector(1,2){5}}
\thicklines\put(32,-21){\vector(1,2){5}}
  \end{picture}}}

%% level 3

 \put(0,-210){\makebox(40,50){
  \begin{picture}(40,50)(0,-50)
\put(0,-2){\makebox(0,0)[t]{1}}
\put(20,-2){\makebox(0,0)[t]{2}}
\put(40,-2){\makebox(0,0)[t]{2}}
\put(10,-22){\makebox(0,0)[t]{1}}
\put(30,-22){\makebox(0,0)[t]{1}}
\put(20,-42){\makebox(0,0)[t]{1}}
\thicklines\put(18,-41){\vector(-1,2){5}}
\thicklines\put(8,-21){\vector(-1,2){5}}
  \end{picture}}}

\put(60,-210){\makebox(40,50){
  \begin{picture}(40,50)(0,-50)
\put(0,-2){\makebox(0,0)[t]{1}}
\put(20,-2){\makebox(0,0)[t]{2}}
\put(40,-2){\makebox(0,0)[t]{2}}
\put(10,-22){\makebox(0,0)[t]{1}}
\put(30,-22){\makebox(0,0)[t]{1}}
\put(20,-42){\makebox(0,0)[t]{2}}
\thicklines\put(22,-41){\vector(1,2){5}}
\thicklines\put(28,-21){\vector(-1,2){5}}
  \end{picture}}}

 \put(140,-210){\makebox(40,50){
  \begin{picture}(40,50)(0,-50)
\put(0,-2){\makebox(0,0)[t]{1}}
\put(20,-2){\makebox(0,0)[t]{2}}
\put(40,-2){\makebox(0,0)[t]{2}}
\put(10,-22){\makebox(0,0)[t]{1}}
\put(30,-22){\makebox(0,0)[t]{2}}
\put(20,-42){\makebox(0,0)[t]{1}}
\thicklines\put(18,-41){\vector(-1,2){5}}
\thicklines\put(8,-21){\vector(-1,2){5}}
  \end{picture}}}

\put(200,-210){\makebox(40,50){
  \begin{picture}(40,50)(0,-50)
\put(0,-2){\makebox(0,0)[t]{1}}
\put(20,-2){\makebox(0,0)[t]{2}}
\put(40,-2){\makebox(0,0)[t]{2}}
\put(10,-22){\makebox(0,0)[t]{1}}
\put(30,-22){\makebox(0,0)[t]{2}}
\put(20,-42){\makebox(0,0)[t]{2}}
\thicklines\put(22,-41){\vector(1,2){5}}
\thicklines\put(32,-21){\vector(1,2){5}}
  \end{picture}}}

 \put(280,-210){\makebox(40,50){
  \begin{picture}(40,50)(0,-50)
\put(0,-2){\makebox(0,0)[t]{1}}
\put(20,-2){\makebox(0,0)[t]{2}}
\put(40,-2){\makebox(0,0)[t]{2}}
\put(10,-22){\makebox(0,0)[t]{2}}
\put(30,-22){\makebox(0,0)[t]{2}}
\put(20,-42){\makebox(0,0)[t]{1}}
\thicklines\put(18,-41){\vector(-1,2){5}}
\thicklines\put(12,-21){\vector(1,2){5}}
  \end{picture}}}

\put(340,-210){\makebox(40,50){
  \begin{picture}(40,50)(0,-50)
\put(0,-2){\makebox(0,0)[t]{1}}
\put(20,-2){\makebox(0,0)[t]{2}}
\put(40,-2){\makebox(0,0)[t]{2}}
\put(10,-22){\makebox(0,0)[t]{2}}
\put(30,-22){\makebox(0,0)[t]{2}}
\put(20,-42){\makebox(0,0)[t]{2}}
\thicklines\put(22,-41){\vector(1,2){5}}
\thicklines\put(32,-21){\vector(1,2){5}}
  \end{picture}}}

%% level 4

 \put(0,-290){\makebox(40,50){
  \begin{picture}(40,50)(0,-50)
\put(0,-2){\makebox(0,0)[t]{2}}
\put(20,-2){\makebox(0,0)[t]{2}}
\put(40,-2){\makebox(0,0)[t]{2}}
\put(10,-22){\makebox(0,0)[t]{1}}
\put(30,-22){\makebox(0,0)[t]{1}}
\put(20,-42){\makebox(0,0)[t]{1}}
\thicklines\put(18,-41){\vector(-1,2){5}}
\thicklines\put(8,-21){\vector(-1,2){5}}
  \end{picture}}}
\put(60,-290){\makebox(40,50){
  \begin{picture}(40,50)(0,-50)
\put(0,-2){\makebox(0,0)[t]{2}}
\put(20,-2){\makebox(0,0)[t]{2}}
\put(40,-2){\makebox(0,0)[t]{2}}
\put(10,-22){\makebox(0,0)[t]{1}}
\put(30,-22){\makebox(0,0)[t]{1}}
\put(20,-42){\makebox(0,0)[t]{2}}
\thicklines\put(22,-41){\vector(1,2){5}}
\thicklines\put(28,-21){\vector(-1,2){5}}
  \end{picture}}}

 \put(140,-290){\makebox(40,50){
  \begin{picture}(40,50)(0,-50)
\put(0,-2){\makebox(0,0)[t]{2}}
\put(20,-2){\makebox(0,0)[t]{2}}
\put(40,-2){\makebox(0,0)[t]{2}}
\put(10,-22){\makebox(0,0)[t]{1}}
\put(30,-22){\makebox(0,0)[t]{2}}
\put(20,-42){\makebox(0,0)[t]{1}}
\thicklines\put(18,-41){\vector(-1,2){5}}
\thicklines\put(8,-21){\vector(-1,2){5}}
  \end{picture}}}
\put(200,-290){\makebox(40,50){
  \begin{picture}(40,50)(0,-50)
\put(0,-2){\makebox(0,0)[t]{2}}
\put(20,-2){\makebox(0,0)[t]{2}}
\put(40,-2){\makebox(0,0)[t]{2}}
\put(10,-22){\makebox(0,0)[t]{1}}
\put(30,-22){\makebox(0,0)[t]{2}}
\put(20,-42){\makebox(0,0)[t]{2}}
\thicklines\put(22,-41){\vector(1,2){5}}
\thicklines\put(32,-21){\vector(1,2){5}}
  \end{picture}}}

 \put(280,-290){\makebox(40,50){
  \begin{picture}(40,50)(0,-50)
\put(0,-2){\makebox(0,0)[t]{2}}
\put(20,-2){\makebox(0,0)[t]{2}}
\put(40,-2){\makebox(0,0)[t]{2}}
\put(10,-22){\makebox(0,0)[t]{2}}
\put(30,-22){\makebox(0,0)[t]{2}}
\put(20,-42){\makebox(0,0)[t]{1}}
\thicklines\put(18,-41){\vector(-1,2){5}}
\thicklines\put(12,-21){\vector(1,2){5}}
  \end{picture}}}
\put(340,-290){\makebox(40,50){
  \begin{picture}(40,50)(0,-50)
\put(0,-2){\makebox(0,0)[t]{2}}
\put(20,-2){\makebox(0,0)[t]{2}}
\put(40,-2){\makebox(0,0)[t]{2}}
\put(10,-22){\makebox(0,0)[t]{2}}
\put(30,-22){\makebox(0,0)[t]{2}}
\put(20,-42){\makebox(0,0)[t]{2}}
\thicklines\put(22,-41){\vector(1,2){5}}
\thicklines\put(32,-21){\vector(1,2){5}}
  \end{picture}

  }}

 \put(180,-310){\makebox(20,10){Fig.7}}
\end{picture}
\end{center}

The diagrams are grouped into four rows with three pairs in each row.
Let us ennumerate the diagrams by triples of numbers $(q_1,q_2,q_3)$, where
$q_3=1,2,3,4$ is the row number, $q_2=1,2,3$ is the
number of the group and $q_1=1,2$ is the number inside the group. For a
diagram with a number $(q_1,q_2,q_3)$ its rows represent the functions
$g_1^{q_1,q_2,q_3}:P_0\to {1,2}, g_2^{q_1,q_2,q_3}:P_1\to {1,2},
g_3^{q_1,q_2,q_3}:P_2\to {1,2}$ (where $P_0, P_1, P_2$ are ordered sets
from Section 3.2 such that $|P_0|=1, |P_1|=2, |P_2|=3$), so to say
``function $g_k^{q_1,q_2,q_3}$'' is the same as to say ``the $k-$th row
of the diagram $(q_1,q_2,q_3)$'' and vice versa. Each triple $(q_1,q_2,q_3)$
is just a sequence $Q$ of $1-$element subsets of $P_1,...,P_3$ (see Section
3.2), so here we can say ``$(q_1,q_2,q_3)$-diagram'' instead of
``$Q-$diagram for $Q=(q_1,q_2,q_3)$''.
The set of $Q-$paths on Fig.7 gives us the diagonal monomial
%(10)
%^%\begin{equation}
$$a_{111}a_{112}a_{121}a_{122}a_{211}a_{212}a_{221}a_{222}
a_{111}a_{122}a_{211}a_{222}a_{111}a_{212}a_{121}a_{222}
a_{111}a_{221}a_{112}a_{222}\times$$
$$\times a^2_{111}a^2_{222}$$
%^%\end{equation}
of $Det(2\times 2\times 2)$. The rest of monomials is obtained as
$Q-$paths over diagrams obtained from initial ones by
permutations of the following type:

i) Permutations of odd type. These are permutations of rows between
different diagrams. The generators are:

1) for a given $(q_2,q_3)$ and $\tau\in S_2$ the action of $\tau$
on $g_1^{\bullet,q_2,q_3}$ (first rows of diagrams $(1,q_2,q_3)$
and $(2,q_2,q_3)$) is the following
%(11)
\begin{equation}
\tau (g_1^{q_1,q_2,q_3})=g_1^{\tau (q_1),q_2,q_3}
\end{equation}

2) for a given $(q_3)$ and $\tau \in S_3$ the action of $\tau$
on $g_1^{\bullet,\bullet,q_3}$ is the following
%(12)
\begin{equation}
\tau (g_2^{\bullet,q_2,q_3})=g_2^{\bullet,\tau(q_2),q_3}
\end{equation}
and permutation $\tau$ for $q_3=2$ has to be chosen the same as for
$q_3=3$. In other words we have to permute the corresponding second
rows of diagrams with $q_3=2$ and $q_3=3$ synchronously, so the permutations
of second rows are the elements of the group $S_3\times S_3\times S_3$.\\
We assign to an odd type permutation the sign wich correspond to its parity.

ii) Permutations of even type. These are permutations of elements in the
rows (elements of $P_k$). The generators are:

for a given $(q_3)$ permute the elements in second rows (elements of sets
$P_1$)
of $(\bullet,2,q_3)-$diagrams (the $\bullet$ means that this permutation does
not depend on the value of $q_1$) so that these permutations for diagrams with
$q_3=2$ and $q_3=3$ coincide (synchronization condition).\\
We assign to even type permutations positive sign.

The permutations of different (odd and even) types commute.
For a given permutation the sign of the monomial computed from $Q-$paths
over $permutation(Q)-$diagrams equals to the sign of odd type component
of the permutation. As a prelude to the general algorithm we can remark
that "synchronization" of permutations on the space of $Q-$diagrams takes
place for the sets of permutations which have the same domaine of values
of $g_3:P_2\to T_3$.

As soon as we can compute the closed determinant for a form of a given
format, its determinant is computed as the quotient of the closed
determinant and the product of all its minors, which are the
determinants of submatricies of smaller formats.

%{\bf Aknoledgements.}
%The author wants to thank Professors
%J.Barros-Neto, Chanillo, P.Landweber, B.Osofsky and E.Taft, whose
%highest professional skills in teaching readily revealed the
%impossibility for him to continue graduate studies at the Department
%of Mathematics of Rutgers University, which gave him time and
%inspiration to finish this work.

\end{document}